\documentclass{gDEA2e}

\usepackage{epstopdf}
\usepackage{subfigure}

\theoremstyle{plain}
\newtheorem{theorem}{Theorem}[section]

\newtheorem{lemma}[theorem]{Lemma}

\theoremstyle{definition}

\newtheorem{example}[theorem]{Example}

\theoremstyle{remark}
\newtheorem{remark}{Remark}
\begin{document}



\title{\textit{Lyapunov direct method for investigating stability of \\
nonstandard finite difference schemes for metapopulation models}}

\author{
\name{Quang A Dang\textsuperscript{a}$^{\ast}$\thanks{$^\ast$Corresponding author. Email: dangquanga@cic.vast.vn}
Manh Tuan Hoang\textsuperscript{b}}
\affil{\textsuperscript{a}Center for Informatics and Computing, Vietnam Academy of Science and Technology\\ (VAST),  18 Hoang Quoc Viet, Cau Giay, Hanoi, Vietnam;\\
\textsuperscript{b}Institute of Information Technology, 
Vietnam Academy of Science and Technology (VAST),
18 Hoang Quoc Viet, Cau Giay, Hanoi, Vietnam}
}
\maketitle

\begin{abstract}
In this paper nonstandard finite difference (NSFD) schemes  of two metapopulation models are constructed. The stability properties of the discrete models are investigated by the use of a generalization of Lyapunov stability theorem. Due to this result we have proved that the NSFD schemes preserve all properties of the metapopulation models. Numerical examples confirm the obtained theoretical results of the properties of the constructed difference schemes. The method of Lyapunov functions proves to be much simpler than the standard method for studying stability of the discrete metapopulation model in our very recent paper. 
\end{abstract}


\begin{classcode}37M05, 39A10, 65L12, 65L20 \end{classcode}

\section{Introduction}
 Many phenomena and processes in physics, chemistry, biology, ecology, finance, environment etc. are modelled by ordinary or partial differential equations
 \cite{Linda, Basak, Barrante, Brauer, Henner, Keshet, Marchuk1, Marchuk2, Robinson, Smith, Stakgold}. The solutions of these equations often possess specific properties, such as positivity, monotonicity, periodicity, stability and some invariant properties. In general, these differential equations are very complicated and it is difficult, even impossible, to find their exact solutions. Therefore, the study of numerical methods and simulation for the solution of these differential equations is very important. Among the numerical methods for differential equations the finite difference method is most popular. The theory of this method for ordinary and partial differential equations is developed rather completed in  \cite{Leveque, Petzold, Samarskii, Strikwerda}. We call the difference schemes in these books and all related sources \emph{standard finite difference schemes}. In many nonlinear problems the standard difference schemes exhibit a serious drawback which is called "numerical instabilities" \cite{Mickens1, Mickens2, Mickens4}. Under this concept we have in mind the phenomena when the discrete models, for example, the difference schemes, do not preserve properties of the corresponding differential equations. In  \cite{Mickens1, Mickens2, Mickens3, Mickens4}  Mickens  showed many examples and analysed the numerical instabilities when using standard difference schemes. In general,  standard difference schemes preserve the properties of the differential equations only in the case if the discretization parameter  $h$ is sufficiently small. Therefore, when studying dynamical models in  large time intervals the selection of small time steps will requires very large computational effort, so these discrete models  are inefficient. Besides, for some special dynamical problems  standard difference schemes cannot preserve the properties of the problems for any step sizes.\par
In order to overcome the numerical instabilities phenomena in 1989  Mickens \cite{Mickens0} introduced the concept \emph{Nonstandard Finite Difference} (NSDF) schemes and after that has developed NSDF methods in many works, such as \cite{Mickens1, Mickens2, Mickens3, Mickens4}. According to Mickens, NSDF schemes are those constructed following a set of five basic rules. The NSDF schemes preserve main properties  of the differential counterparts, such as positivity, monotonicity, periodicity, stability and some other invariants including energy and geometrical shapes. It should be emphasized that NSFD schemes can preserve all properties of the continuous models for any discretization parameters. The discrete models with these properties are called \emph{dynamically consistent}  \cite{DQA, AL2, DK4, Garba, Mickens5, Roeger5, Roeger6}.\par

For the last two decades NSDF methods have attracted attention from many researchers and achieved significant results  \cite{AL1, DK3,DK5, E1, K1, K2, Partidar, Roeger1, Roeger2, Roeger4, Roeger5, Wood}. The property of stability of the set of equilibria of differential equations is one of these results because it plays the essential role in the study of asymptotical behaviour of the solutions of differential equations. The construction of difference schemes, which preserve the stability of the equilibrium points, is important in numerical simulation of differential equations. The difference schemes with this stability property is called \emph{elementary stable} schemes  \cite{AL1, DK1, DK2, Wood}. There are many works concerning the elementary stable schemes. The typical results are for general dynamical systems \cite{DK1, DK2} and for other specific systems \cite{DK5, Roeger3, Wood} etc. One popular approach to the elementary stability is  the investigation of Jacobian matrice of the discrete models at the equilibria, namely, determination of conditions ensuring that all eigenvalues of Jacobian matrice have moduli less or equal to 1. This is the necessary and sufficient condition for hyperbolic equilibrium points to be locally stable \cite{Linda, Keshet}. The mentioned above approach has the following weaknesses and limitations:\par
\begin{enumerate}
\item It is applicable when all the equilibrium points are hyperbolic. To our best knowledge, at present no results on NSFD schemes preserving the stability of non-hyperbolic equilibrium points are available.
\item Even when all the equilibrium points are hyperbolic, the determination of the conditions ensuring all the eigenvalues of Jacobian to be in the unit ball, is very difficult. Theoretically, it is possible to use the Jury's conditions \cite{Linda, Elaydi, Keshet} for finding these conditions, but in practice, this is extremely complex in many cases, for example, when  the system of equations has large dimension or contains several parameters.
\item The consideration of Jacobian only guarantees the local stability meanwhile  many models have the global stability.
\end{enumerate}

In order to overcome the above shortcomings of the approach with the use of Jacobian it is useful to use the Lyapunov stability theory for proving the stability of both hyperbolic and non-hyperbolic equilibrium points  \cite{Linda, Elaydi}. The drawback of this approach is that not always it is possible to find appropriate Lyapunov's function. Nevertheless, in many specific problems it is easy to find associated Lyapunov's functions. Then the stability of the equilibrium points may be established without the study of Jacobian matrices of discrete models. This is a perspective approach for many discrete models.\par

For illustrating this approach in this paper we consider two metapopulation models, one was proposed by Keymer \cite{Keymer} and another by Amarasekare \cite{Amarasekare}. This is the models with complex properties. By using a generalization of Lyapunov's stability theorem \cite{Iggidr} we construct NSDF schemes preserving stability properties of the models. This way is much simpler than the using Jacobian because it doesn't require complicated calculations and difficult techniques.\par

The paper is organized as follows. In Section 2 we recall the two models of metapopulations. In Sections 3 and 4 we construct NSDF schemes dynamically consistent  with the continuous models. Next, in Section 5 we report some numerical experiments for validating the obtained theoretical results. Finally, some concluding remarks are given in Section 6.
\section{Mathematical models of metapopulations}
\subsection{Keymer's metapopulation model }
Consider the metapopulation model proposed by Keymer in 2000  \cite{Keymer}. It is described by the system of three nonlinear differential equations
\begin{equation}\label{eq:1}
\dfrac{dp_0}{dt} = e(p_1 + p_2) - \lambda p_0, \quad \dfrac{dp_1}{dt} = \lambda p_0 - \beta p_1 p_2 + \delta p_2 - ep_1, \quad \dfrac{dp_2}{dt}= \beta p_1 p_2 - (\delta + e)p_2, 
\end{equation}
where $p_0$, $p_1$ and $p_2$ denote the proportion of uninhabitable patches, the proportion of the habitable patches that are not occupied and the proportion of habitable patches that are occupied, respectively, $\lambda$ is the rate of patch creation, $e$ is the rate of patch destruction, $\delta$ is the rate of population extinction and $\beta$ is the rate of propagule reproduction. Because of $p_0 + p_1 + p_2 = 1$, the system \eqref{eq:1} can be reduced to the two equations
\begin{equation}\label{eq:2}
\dfrac{dp_1}{dt} = \lambda(1 - p_1 - p_2) - \beta p_1 p_2 + \delta p_2 - ep_1, \qquad \dfrac{dp_2}{dt} = p_2(\beta p_1 - \delta - e). 
\end{equation}
From the biological meaning of the model we shall consider the initial conditions $p_1(0), p_2(0)$   satisfying
\begin{equation}\label{eq:3}
\big(p_1(0), p_2(0)\big) \in D_2 := \Big\{(x, y) \in \mathbb{R}^2: 0 \leq x, y; x + y \leq 1\Big\}.
\end{equation}
The mathematical analysis in  \cite{Linda, Keymer} shows that the model \eqref{eq:2} possesses the following properties:
\begin{description}
\item[$(P_1)$] The monotone convergence of the sum $s(t) := p_1(t) + p_2(t)$:\\
With the initial conditions satisfying \eqref{eq:3} the sum of the solutions $s(t) := p_1(t) + p_2(t)$ monotonically converges to $s^* := \lambda/(\lambda + e)$.
\item[$(P_2)$] Boundedness\\ All the solution $p_1(t), p_2(t)$ with the initial conditions satisfying \eqref{eq:3} also satisfy \eqref{eq:3}. In other words, the set $D_2$ is positive invariant.
\item[$(P_3)$] Local asymptotic stability\\The model \eqref{eq:2} has two equilibria
\begin{equation*}
P_1^* = \big(\dfrac{\lambda}{\lambda + e}, 0\big), \qquad P_2^* = (\dfrac{\delta + e}{\beta}, 1 - \dfrac{e}{\lambda + e} - \dfrac{\delta + e}{\beta}) = (\dfrac{\delta + e}{\beta}, \dfrac{\lambda}{\lambda + e} - \dfrac{\delta + e}{\beta}).
\end{equation*}
Set $\mathcal{R}_0 = \dfrac{\beta \lambda}{(\lambda + e)(\delta + e)}$. This is a threshold for persistence indicating the number of propagules needed during the species and the patch lifetime for the species to persist for the model.

According to \cite[Chapter 6]{Linda} the first equilibrium $P_1^*$ is locally asymptotically stable if  $\mathcal{R}_0 < 1$ and the second equilibrium $P_2^*$ is locally asymptotically stable if $\mathcal{R}_0 > 1$.
\item[$(P_4)$] Global stability\\ Due to the Poincare - Bendixson theory \cite[Chapter 5-6]{Linda} it is possible to show that if $\mathcal{R}_0 > 1$ then the second equilibrium is globally asymptotically stable, and if $\mathcal{R}_0 < 1$ then the first equilibrium is globally asymptotically stable. 

\item[$(P_5)$] Non-periodic solution\\ Applying Dulax's criterion it is possible to show that the system \eqref{eq:2} does not have periodic solutions in the domain $D = \{0 < p_1 + p_2 < 1\}$.
\end{description}
Clearly, this is a model with complex properties. Very recently, in \cite{DQA} by using NSFD methods we have successfully constructed a discrete metapopulation model dynamically consistent with the continuous counterpart. It means that the discrete model preserves all the properties 
$(P_1)-(P_5)$ of the model \eqref{eq:1}. The construction of the NSFD scheme is very complicated. The main difficulty is in the stability properties  $(P_3)$ and $(P_4)$. For constructing NSFD scheme preserving these properties it is needed to use difficult techniques and cumbersome calculations. \par
  In Section 3 by using the Lyapunov's stability theory we construct NSFD scheme preserving the global stability property of the model \eqref{eq:1}, so due to it the difference scheme preserves all the properties of the model. As will see later, the method of  Lyapunov's function is much simpler and shorter than the method for proving stabilities in   \eqref{eq:1}.

\subsection{Amarasekare-Possingham's metapopulation model}
Consider the metapopulation model proposed by Amarasekare and Possingham in 2001 \cite{Amarasekare}. It is described by the system of four nonlinear differential equations
\begin{equation}\label{eq:4}
\begin{split}
&\dfrac{dI}{dt} = \beta_I SI - e_I I + fL - gI, \qquad \qquad \dfrac{dS}{dt} = e_I I - \beta_I SI + fR - gS,\\
&\dfrac{dL}{dt} = gI - fL - e_LL + \beta_LRI, \qquad \qquad \dfrac{dR}{dt} = gS - fR + e_LL - \beta_LRI.\\
\end{split}
\end{equation}

Here $f$ is the disturbance frequency and $g$, the rate of habitat succession. Quantities $e_I$ and $e_L$ represent local extinction rates, and $\beta_I$ and $\beta_L$ the per patch colonization rates of infected and latent patches, respectively. The total number of patches in the system is assumed to be constant such that $I + S + L + R = P$. Alternatively, $I, S, L$ and $R$ can be thought of as the frequency of each patch type in the landscape in which case $I + S + L + R = 1$.\par

From the biological meaning of the model we shall consider the initial conditions $I(0), S(0)$, $L(0), R(0)$   satisfying
\begin{equation}\label{eq:5}
\big(I(0), S(0), L(0), R(0)\big) \in D_4 := \big\{(I, S, L, R) \in \mathbb{R}^4: 0 \leq I, S, L, R; I + S + L + R = 1\Big\}.
\end{equation}
The mathematical analysis shows that the model \eqref{eq:4} possesses the following properties:
\begin{description}
\item[$(P_1)$] The monotone convergence of the sum $a(t) := I(t) + S(t)$ and $b(t) := L(t) + R(t)$:\\
For any initial conditions satisfying
 \eqref{eq:5} the sum $a(t) := I(t) + S(t)$  monotonically converges to $a^* := f/(f + g)$, and $b(t) := L(t) + R(t)$  monotonically converges to $b^* := g/(f + g)$.
\item[$(P_2)$] Boundedness:\\ 
All the solutions $I(t), S(t), L(t), R(t)$ with the initial conditions satisfying \eqref{eq:5} also satisfy \eqref{eq:5}. In other words, the set $D_4$ is positive invariant. 
\item[$(P_3)$] Local asymptotic stability (see \cite{Amarasekare}):\\
The model \eqref{eq:5} has two equilibria (lying on the boundary or inside of $D_4$) $E_i^* = (I_i^*, S_i^*, L_i^*, R_i^*), i = 1, 2$, where
\begin{equation}\label{eq:6}
I_1^* = 0, \qquad S_1^* = \dfrac{f}{f + g}, \qquad L_1^* = 0, \qquad R_1^* = \dfrac{g}{g + f}.
\end{equation}
\begin{equation}\label{eq:7}
\begin{split}
&I_2^* := I^* = \dfrac{-b + \sqrt{b^2 - 4ac}}{2a}, \quad S_2^* := S^* = \dfrac{b + 2\beta_I\beta_L \dfrac{f}{f + g} - \sqrt{b^2 - 4ac}}{2a}, \\
&R_2^* := R^* = \dfrac{g}{f + g} - \dfrac{\beta_I}{f}{I^*}^2 + \big(\dfrac{\beta_I}{f + g} - \dfrac{g + e_I}{f}\big)I^*,\\
&L_2^* := L^* = 1 - I^* - S^* - R^* = \dfrac{g}{f + g} - R^* = \dfrac{\beta_I}{f}{I^*}^2 - \big(\dfrac{\beta_I}{f + g} - \dfrac{g + e_I}{f}\big)I^*,
\end{split}
\end{equation}
\begin{equation}\label{eq:7a}
\begin{split}
a =& \beta_I \beta_L, \qquad b = \beta_I(f + e_L) + \beta_L(e_I + g) - \beta_I\beta_L\dfrac{f}{f + g},\\
c =& (f + e_L)(e_I - \beta_I\dfrac{f}{f + g}) + g(e_L - \beta_L \dfrac{f}{f + g}).
\end{split}
\end{equation}
Set  $\mathcal{R}_0 =  1 - c$. Then  $\mathcal{R}_0$ is threshold parameter, i.e., if $\mathcal{R}_0 > 1$ then the model has   positive equilibrium point, otherwise positive equilibrium point does not exist. The first equilibrium point is locally asymptotically if and only if  $c > 0 \,\, (\mathcal{R}_0 < 1)$ and the second equilibrium point is locally asymptotically if and only if $c < 0 \,\, (\mathcal{R}_0 > 1)$ .
\end{description}
\section{NSFD scheme for the model \eqref{eq:2}}
In this section we construct NSFD scheme for the model \eqref{eq:2} so that the obtained difference scheme preserves all dynamic properties of the original continuous model for any discretization parameter or step size $h > 0$. According to Mickens, a finite difference scheme is called \emph{nonstandard} if at least one of the following conditions is satisfied 
 \cite{Mickens1, Mickens2, Mickens3, Mickens4}:
\begin{itemize}
\item A nonlocal approximation is used.
\item  The discretization of the derivative is not traditional and uses a function $0 < \varphi(h) = h + \mathcal{O}(h^2)$.
\end{itemize} 

For simplicity of presentation, in the model \eqref{eq:2} we use the notations $x(t)$ and $y(t)$ instead of  $p_1(t)$ and $p_2(t)$, respectively. We reconsider the family of difference schemes of the form  \cite{DQA}

\begin{equation}\label{eq:8}
\begin{split}
\dfrac{x_{k + 1} - x_k}{\varphi(h)} &= -c_1(\lambda + e)x_k - c_2(\lambda + e)x_{k + 1} + c_3(\delta - \lambda)y_k + c_4(\delta - \lambda)y_{k + 1}\\
& - c_5\beta x_ky_k - c_6\beta x_{k + 1}y_k - c_7\beta x_ky_{k + 1} - c_8 \beta x_{k + 1}y_{k + 1}+ \lambda,\\
\dfrac{y_{k + 1} - y_k}{\varphi(h)} & = -c_1(\lambda + e)y_k - c_2(\lambda + e)y_{k + 1} + c_3(\lambda - \delta)y_k + c_4(\lambda - \delta)y_{k + 1}\\
 &+ c_5\beta x_ky_k + c_6\beta x_{k + 1}y_k + c_7\beta x_ky_{k + 1} + c_8\beta x_{k + 1}y_{k + 1},
\end{split}
\end{equation}
where the parameters $c_i$ and the function $\varphi(h)$ satisfy
\begin{equation}\label{eq:9}
 c_1 + c_2 = 1, \quad c_3 = 1, \quad c_4 = 0, \quad c_5 + c_6 = 1, \quad c_7 = c_8 = 0, \quad \varphi(h) = h + \mathcal{O}(h^2).
\end{equation}
The explicit form of the difference scheme is defined by \cite[Theorem 2.4]{DQA}
\begin{equation}\label{eq:1d}
\begin{split}
&x_{k + 1} = \dfrac{x_k - \varphi c_1 (\lambda + e)x_k + \varphi (\delta - \lambda)y_k - \varphi c_5 \beta x_ky_k + \lambda \varphi}{1 + \varphi c_2(\lambda + e) + \varphi c_6 \beta y_k},\\
&y_{k + 1} =  \dfrac{[1 + \varphi c_2 (\lambda + e) + \varphi c_6 \beta y_k][y_k - \varphi c_1 (\lambda + e)y_k + \varphi (\lambda - \delta)y_k + \varphi c_5 \beta x_ky_k]}{[1 + \varphi c_2(\lambda + e) + \varphi c_6 \beta y_k][1 + c_2(\lambda + e)]}\\
 &+ \dfrac{\varphi c_6 \beta y_k[x_k - \varphi c_1 (\lambda + e)x_k + \varphi (\delta - \lambda)y_k - \varphi c_5 \beta x_ky_k + \lambda \varphi]}{[1 + \varphi c_2(\lambda + e) + \varphi c_6 \beta y_k][1 + c_2(\lambda + e)]}.
\end{split}
\end{equation}

\subsection{Properties $(P_1)$ and $(P_2)$}
The results of the difference schemes preserving the properties 
 $(P_1)$ and $(P_2)$ of the model \eqref{eq:2} 
 are stated in the following theorems   \cite{DQA}:
\begin{theorem}\label{Theorem1}
The difference scheme \eqref{eq:8}-\eqref{eq:9} preserves Property $(P_1)$ of the model \eqref{eq:2} if
\begin{equation}\label{eq:10}
c_1 \leq 0, \qquad c_2 \geq 0.
\end{equation}
\end{theorem}
\begin{theorem}\label{Theorem2}
Consider the scheme \eqref{eq:8}-\eqref{eq:9}. Under the assumptions
\begin{equation}\label{eq:11}
c_5 \leq 0, \quad c_6 \geq 0, \quad c_2 \geq c_6, \quad c_1 \leq -\dfrac{\delta}{\lambda + e},
\end{equation}
the scheme \eqref{eq:8}-\eqref{eq:9} preserves Property $(P_2)$ of the model \eqref{eq:2}.
\end{theorem}

\subsection{Stability properties $(P_3)$ and $(P_4)$}
The main difficulty in the construction of difference schemes preserving the properties of the model \eqref{eq:2} is in Properties $(P_3)$ and $(P_4)$, especially in the property of global stability $(P_4)$. In this section by using a generalization of Lyapunov's stability theorem  \cite[Theorem 3.3]{Iggidr} we show the stability properties of the equilibrium points in a simple and easy way.

For this purpose we consider the function
\begin{equation}\label{eq:11}
V(x, y) = \big(x + y - \dfrac{\lambda}{\lambda + e}\big)^2, \qquad (x, y) \in D_2,
\end{equation}
where $D_2$ is defined by \eqref{eq:3}.
We shall show that the function
 $V(x, y)$ satisfies all the conditions of \cite[Theorem 3.3]{Iggidr} on $D_2$. Obviously, $V(x, y)$ is continuous on $D_2$. Moreover,\\
(i) $V(x, y) \geq 0$ for any $(x, y) \in D_2$ and  $V(P_i^*) = 0$ $(i = 1, 2)$.\\
(ii) \begin{equation*}
\begin{split}
\Delta V(x_k, y_k) &= V(x_{k + 1}, y_{k + 1}) - V(x_k, y_k) \\
&= \big(x_{k + 1} + y_{k +1} - x_k - y_k\big)\big(x_{k + 1} + y_{k + 1} + x_k + y_k - 2\dfrac{\lambda}{\lambda + e} ).
\end{split}
\end{equation*}
Since Property 
 $(P_1)$ of the model \eqref{eq:2} is preserved it follows that $\Delta V(x_k, y_k) \leq 0$ for any $(x_k, y_k) \in D_2$.  Besides $\Delta V(x_k, y_k) = 0$ if and only if  $x_k + y_k = \dfrac{\lambda}{\lambda + e}$. Therefore,
 \begin{equation*}
G^* = \big\{(x, y) \in D_2: x + y = \dfrac{\lambda}{\lambda + e}\big\},
\end{equation*}
where $G^*$ is the largest positively invariant set containing in $G=\{(x,y): \Delta V(x,y)=0  \}$.\\
(iii) On the other hand, due to the preservation of Property  $(P_2)$ of the model all the solutions of \eqref{eq:8} are bounded. Hence, all the conditions (1), (2) and (4) of \cite[Theorem 3.3]{Iggidr} are satisfied. In order to prove the global stability of \eqref{eq:2} it remains to show that
\begin{enumerate}
\item If $\mathcal{R}_0 < 1$ then $P_1^*$ is $G^*-$ globally asymptotically stable, i.e., $P_1^*$ is $G^*$- asymptotically stable and $G^*$- globally attractive.
\item If $\mathcal{R}_0 > 1$ then $P_2^*$ is $G^*-$ globally asymptotically stable, i.e.,  $P_2^*$ is $G^*$- asymptotically stable and $G^*$- globally attractive.
\end{enumerate}
Since Property $(P_1)$ of the model \eqref{eq:2} is preserved for any initial conditions belonging to $G^*$, i.e., $x_0 + y_0 = \dfrac{\lambda}{\lambda + e}$, the solutions of \eqref{eq:8}-\eqref{eq:9} satisfy $x_k + y_k = \lambda/(\lambda + e)$. Substituting $y_k = \lambda/(\lambda + e) - x_k$ into \eqref{eq:1d} we obtain the scheme depending only on  $x_k$ of the form
\begin{equation}\label{eq:12}
x_{k + 1} = x_k + \dfrac{\varphi \beta x_k^2 - \varphi \Big[\dfrac{\beta \lambda}{\lambda + e} + \big(\delta + e\big)\Big]x_k + \varphi\lambda \dfrac{\delta + e}{\lambda + e}}{1 + \varphi c_2(\lambda + e) + \varphi \beta c_6 \big(\dfrac{\lambda}{\lambda + e} - x_k\big)}.
\end{equation}
This is the discretization of the model \eqref{eq:2} on the set $\Big\{\big(x(t), y(t)\big) \in D_2: x(t) + y(t) = \lambda/(\lambda + e)\Big\}$. Denote by
\begin{equation}\label{eq:1a}
P_{1, x}^* = \dfrac{\lambda}{\lambda + e}, \qquad P_{2, x}^* = \dfrac{\delta + e}{\beta},
\end{equation}
 $x-$components of  the equilibrium points $P_i^* \ (i = 1, 2)$. This is also two equilibrium points of \eqref{eq:12}. Therefore, the proof of $G*$- global stability of the equilibrium points $P_i^* \ (i = 1, 2)$ is equivalent to showing that
 \begin{enumerate}
\item If $\mathcal{R}_0 < 1$ then $P_{1, x}^*$ is the globally asymptotically stable equilibrium point of  \eqref{eq:12} in $D_2^* := \Big\{x \in \mathbb{R}: 0 \leq x \leq \dfrac{\lambda}{\lambda + e}\Big\}$.
\item If $\mathcal{R}_0 > 1$ then $P_{2, x}^*$ is the globally asymptotically stable equilibrium point of   \eqref{eq:12} in $D_2^*$.
\end{enumerate}

Of course, it is possible to show the local stability of $P_{1, x}^*$ and $P_{2, x}^*$ via the Jacobian matrice of the scheme  \eqref{eq:12}, and doing this is not difficult. However, the scheme \eqref{eq:12} depends only on  $x_k$. Therefore, it is wise to use directly the results of 
 Lubuma  Anguelov \cite[Theorem 3]{AL1} to show the elementary stability of \eqref{eq:12}. 
 Since the set of equilibrium points is preserved it suffices to show the monotone dependence on initial value of the scheme \eqref{eq:12}.
\begin{theorem}\label{Theorem3}
Consider the scheme \eqref{eq:8}-\eqref{eq:9} under the conditions of Theorem \ref{Theorem2}. If additionally assume that
\begin{equation}\label{eq:1b}
c_2 > c^* : = \dfrac{\dfrac{\beta \lambda}{\lambda + e} + \delta + e + c_6\beta(2 - \dfrac{\lambda}{\lambda + e})}{\lambda + e},
\end{equation}
then the scheme \eqref{eq:12} is elementary stable, i.e., $P_{1, x}^*$ and $P_{2, x}^*$ are locally asymptotically stable. 
Hence, the local stability Property $(P_3)$ of the model \eqref{eq:2} is preserved. 
\end{theorem}
\begin{proof}
Consider the function
\begin{equation}\label{eq:1q}
f(x, h) := x + \dfrac{\varphi \beta x^2 - \varphi \Big[\dfrac{\beta \lambda}{\lambda + e} + \big(\delta + e\big)\Big]x + \varphi\lambda \dfrac{\delta + e}{\lambda + e}}{1 + \varphi c_2(\lambda + e) + \varphi \beta c_6 \big(\dfrac{\lambda}{\lambda + e} - x\big)}.
\end{equation}
It is easy to obtain
\begin{equation}\label{eq:1c}
\begin{split}
&\dfrac{\partial f(x, h)}{\partial x} = \dfrac{\big[1 + \varphi c_2 (\lambda + e) + \varphi \beta c_6 \dfrac{\lambda}{\lambda + e}\big] \varphi \big[c_2(\lambda + e) + \beta c_6 \dfrac{\lambda}{\lambda + e} - 2\beta c_6 x - (\dfrac{\beta \lambda}{\lambda + e} + \delta + e)\big]}{\big[1 + \varphi c_2(\lambda + e) + \varphi \beta c_6 \dfrac{\lambda}{\lambda + e} - \varphi c_6 \beta x\big]^2}\\
&+ \dfrac{\big[1 + \varphi c_2 (\lambda + e) + \varphi \beta c_6 \dfrac{\lambda}{\lambda + e}\big](1 + 2 \varphi \beta x)}{\big[1 + \varphi c_2(\lambda + e) + \varphi \beta c_6 \dfrac{\lambda}{\lambda + e} - \varphi c_6 \beta x\big]^2}+ \dfrac{\varphi ^2 \beta^2 (c_6^2 - c_6)x^2 + \varphi^2c_6\beta\lambda \dfrac{\delta + e}{\lambda + e}}{\big[1 + \varphi c_2(\lambda + e) + \varphi \beta c_6 \dfrac{\lambda}{\lambda + e} - \varphi c_6 \beta x\big]^2}.
\end{split}
\end{equation}
Since $c_6 > 1$ there holds $c_6^2 > c_6$. On the other hand, in view of $0 \leq x \leq 1$ from \eqref{eq:1c} and \eqref{eq:1b} it follows
\begin{equation*}
\dfrac{\partial f(x, h)}{\partial x} > \dfrac{\big[1 + \varphi c_2 (\lambda + e) + \varphi \beta c_6 \dfrac{\lambda}{\lambda + e}\big] \varphi (\lambda + e)(c_2 - c^*)}{\big[1 + \varphi c_2(\lambda + e) + \varphi \beta c_6 \dfrac{\lambda}{\lambda + e} - \varphi c_6 \beta x\big]^2} > 0.
\end{equation*}
Therefore, \eqref{eq:12} is monotonically depends on the initial value. The proof of the theorem is complete by using Theorem \cite[Theorem 3]{AL1} .
\end{proof}
The global attractiveness is obtained with the use of following theorem.
\begin{theorem}\label{Theorem4}
Consider the scheme \eqref{eq:8}-\eqref{eq:9}. Under the assumptions of Theorems \ref{Theorem2} and \ref{Theorem3} then
\begin{enumerate}
\item If $\mathcal{R}_0 < 1$ then $P_{1, x}^*$  is globally stable equilibrium point of \eqref{eq:12} in $D_2^*$.
\item If $\mathcal{R}_0 > 1$ then $P_{2, x}^*$ is globally stable equilibrium point of \eqref{eq:12} in $D_2^*$.
\end{enumerate}
\end{theorem}

\begin{proof}
(i) \textbf{Case $\mathcal{R}_0 < 1$}.\\
If $\mathcal{R}_0 < 1$ then $P_{2, x}^* = \dfrac{\delta + e}{\beta} > \dfrac{\lambda}{\lambda + e}$. Therefore $P_{2, x}^* = \dfrac{\delta + e}{\beta}$ 
cannot be a globally stable point of  \eqref{eq:12} in $D_2^*$. We shall show that $P_{1, x}^*$ is a globally stable point of  \eqref{eq:12}.\\
Set $u_k = {\lambda}/{(\lambda + e)} - x_k$. Noticing that on $D_2^*$ there holds $u_k \geq 0$ for any $k$, from \eqref{eq:12} we obtain 
\begin{equation*}
u_{k + 1} = u_k - \dfrac{u_k^2 - \varphi \big[\dfrac{\beta \lambda}{\lambda + e} - (\delta + e)\big]u_k}{1 + \varphi c_2(\lambda + e) + \varphi c_6\beta u_k}.
\end{equation*}
Since $\mathcal{R}_0 < 1$ there holds $\dfrac{\beta \lambda}{\lambda + e} - (\delta + e) < 0$. It follows  $u_{k + 1} < u_k$ for any $k$. The sequence $\big\{u_k\big\}$ is decreasing and bounded from below, consequently, it is convergent. Therefore, the sequence $\big\{x_k\big\}$ also is convergent. Since $\big\{x_k\big\}$ cannot converge to $P_{2, x}^* = \dfrac{\delta + e}{\beta}$ it implies that $P_{1, x}^* = \dfrac{\lambda}{\lambda + e}$ is globally attractive point of \eqref{eq:12} in $D_2^*$.\\
(ii) \textbf{Case $\mathcal{R}_0 > 1$}.\\
For simplicity, here instead of $f(x,h)$ we write $f(x)$.
First, from \eqref{eq:1q} we have
\begin{equation*}
f(x) > x \Longleftrightarrow g(x) := \varphi \beta x^2 - \varphi \Big[\dfrac{\beta \lambda}{\lambda + e} + \big(\delta + e\big)\Big]x + \varphi\lambda \dfrac{\delta + e}{\lambda + e} > 0.
\end{equation*}
The function $g(x)$ has two positive roots $P_{2, x}^* < P_{1, x}^*$ (because $\mathcal{R}_0 > 1$). The equation $g'(x) = 0$ has a unique positive root $x^* = \dfrac{P_{1, x}^* + P_{2, x}^*}{2}$. From the investigating the behaviour of $g(x)$ we obtain
\begin{equation*}
\begin{cases}
f(x) > x \Longleftrightarrow 0 < x < P_{2, x}^*,\\
\\
f(x) < x \Longleftrightarrow  x > P_{2, x}^*.\\
\end{cases}
\end{equation*}
Besides, by Theorem \ref{Theorem3} the function $f(x)$ is increasing. Therefore, from the result \cite[Problem 2.5.38]{Kaczor} it follows that $P_{2, x}$ is a globally attractive point of  \eqref{eq:12}. The proof is complete.
\end{proof}

\begin{remark}{}
The function $V(x,y)$ constructed in the proof of the stability properties of the set of equilibrium points does not satisfy the conditions of the classical Lyapunov stability theorem. It satisfies a generalization of this theorem, namely  \cite[Theorem 3.3]{Iggidr}. The properties $(P_1), (P_2)$ play  an important role in the construction of the function $V(x,y)$ for ensuring the stability of the NSFD scheme.
\end{remark}
\subsection{Non-periodicity of solution}
Non-periodicity of solution of  $(P_5)$ follows from the proved fact that the equilibrium points are globally asymptotically stable. \par
Summarizing the above results we obtain the difference schemes preserving the properties of the model \eqref{eq:2}.
\begin{theorem}\label{Theorem5}
The difference schemes \eqref{eq:8}-\eqref{eq:9} preserve Properties $(P_1)-(P_5)$  of the model \eqref{eq:2} if
\begin{equation}\label{eq:14}
c_5 \leq 0, \quad c_6 \geq 0, \quad c_2 \geq \max\big\{c_6, c^*\big\}, \quad c_1 \leq -\dfrac{\delta}{\lambda + e}.
\end{equation}
\end{theorem}
 \begin{remark}
Theorems 2.22, 2.23, 2.24 in \cite{DQA} and Theorem 5 above give sufficient conditions for the NSFD schemes to preserve the properties of the metapopulation model \eqref{eq:2}. Similar to the scheme in Theorem 2.22 \cite{DQA}, the scheme in Theorem 5 contains 4 parameters $c_1, c_2, c_5, c_6$. The conditions of the both theorems are satisfied for  sufficiently small negative numbers $c_1, c_5$ and  sufficiently large positive numbers $c_2, c_6$, and in this case the constructed difference schemes coincide. Nevertheless, Theorem 2.22 \cite{DQA} only guarantees that the scheme preserves Properties $(P_1)-(P_3)$ of the model although numerical experiments show that Properties $(P_4)$ and $(P_5)$ also are preserved. Moreover, the derivation of the conditions in Theorem 2.22 is very hard because it is done via the analysis of the Jacobian for proving the local stability property $(P_3)$. Using the approach in \cite{DQA} it is difficult to find the conditions for the scheme to preserve Property $(P_4)$ of the model.\par
The NSFD schemes satisfying Theorems 2.23, 2.24 in \cite{DQA} and Theorem 5 preserve Properties $(P_1)-(P_5)$ of the model. Due to the complexity of the model the scheme in Theorem 2.23 contains only two parameters $c_1$ and $c_2$, and the scheme in Theorem 2.24 contains only one parameter being the denominator function $\varphi$. The scheme satisfying Theorem 5 does not satisfy Theorem 2.23 and Theorem 2.24, and conversely. However, the scheme in Theorem 5 still contains 4 parameters, and this allows to construct schemes with  other properties, for example, the second order accuracy property. Besides, it should be emphasized that the derivation of conditions in Theorem 2.23 and Theorem 2.24 is very complicated, while the conditions of Theorem 5 are easily derived by using a generalization of Lyapunov stability theorem.
 \end{remark}
\section{NSFD schemes for the model \eqref{eq:4}}
We propose NSFD schemes for the model \eqref{eq:4} in the form
\begin{equation}\label{eq:15}
\begin{split}
\dfrac{S_{k + 1} - S_k}{\varphi} & = e_I I_k - \beta_I S_{k + 1}I_k + fR_k - gS_{k},\\
\dfrac{I_{k + 1} - I_k}{\varphi} & = \beta_I S_{k + 1}I_k - e_I I_k + fL_K - gI_{k },\\
\dfrac{R_{k + 1} - R_k}{\varphi} & = gS_{k } - fR_k + e_LL_k - \beta_L R_{k + 1}I_k,\\
\dfrac{L_{k + 1} - L_k}{\varphi} &= gI_{k } - fL_k - e_LL_{k } + \beta_L R_{k + 1}I_k.
\end{split}
\end{equation}
Our task now is to determine the conditions for the function $\varphi(h)$ so that the scheme \eqref{eq:15} preserves Properties  $(P_1) - (P_3)$ of the model \eqref{eq:4}.
\subsection{Monotone convergence}
\begin{theorem}\label{theorem5}
The scheme \eqref{eq:15} preserves Property $(P_1)$ of the model \eqref{eq:4} if the function $\varphi(h)$ satisfies
\begin{equation}\label{eq:16a}
\varphi(h) < \dfrac{1}{f + g}, \qquad \forall h > 0.
\end{equation}
\end{theorem}

\begin{proof}
Set $a_k := S_k + I_k$, $b_k := R_k + L_k$. In \eqref{eq:15} adding consecutively the first equation with the second one, the third equation with the fourth one we obtain
\begin{equation*}
a_{k + 1} = (1 - \varphi f - \varphi g)a_k + \varphi f, \qquad b_{k + 1} = (1 - \varphi f - \varphi g)b_k + \varphi g.
\end{equation*}
From here we have
\begin{equation}\label{eq:17}
a_k = \Big(a_0 - \dfrac{f}{f + g}\Big)\Big(1 - \varphi f - \varphi g\Big)^k + \dfrac{f}{f + g},\quad
b_k =  \Big(b_0 - \dfrac{g}{f + g}\Big)\Big(1 - \varphi f - \varphi g\Big)^k + \dfrac{g}{f + g}.
\end{equation}
Since $\varphi(h)$ satisfies \eqref{eq:16a} then $ 1 - \varphi f - \varphi g \in (0, 1)$. From here it follows the proof of the theorem.
\end{proof}

\subsection{Boundedness}
\begin{theorem}\label{theorem5}
The scheme \eqref{eq:15} preserves Property $(P_2)$ of the model \eqref{eq:4} if the function $\varphi(h)$ satisfies
\begin{equation}\label{eq:16}
\varphi(h) < \min\Bigg\{\dfrac{1}{e_I + g}, \quad \dfrac{1}{f + e_L}\Bigg\}, \qquad \forall h > 0.
\end{equation}
\end{theorem}
\begin{proof}
We prove the theorem by induction. First, adding the equations of  \eqref{eq:15} side-by-side we obtain
\begin{equation*}
S_{k + 1} + I_{k + 1} + R_{k + 1} + L_{k + 1} = S_k + I_k + R_k + L_k.
\end{equation*}
Therefore, if $S_k + I_k + R_k + L_k = 1$ then 
\begin{equation}\label{eq:16new}
S_{k + 1} + I_{k + 1} + R_{k + 1} + L_{k + 1} = 1.
\end{equation}
On the other hand, it is easy to transform the scheme \eqref{eq:15} to the explicit form
\begin{equation}\label{eq:17}
\begin{split}
&S_{k + 1} = \dfrac{(1 - \varphi g)S_k + \varphi e_I I_k + \varphi f R_k}{1 + \varphi \beta_I I_k}, \qquad I_{k + 1} = {(1 - \varphi e_I - \varphi g)I_k + \varphi \beta_I S_{k + 1}I_k + \varphi f L_k},
\\
&R_{k + 1} = \dfrac{(1 - \varphi f)R_k +\varphi gS_{k} + \varphi e_LL_k}{1 + \varphi \beta_L I_k}, \qquad L_{k + 1} = {(1 - \varphi f - \varphi e_L)L_k + \varphi g I_{k} + \varphi \beta_L R_{k + 1}I_k}.
\end{split}
\end{equation}
Hence, if $S_k, I_k, R_k, L_k \geq 0$ and $\varphi(h)$ satisfy \eqref{eq:16} then $S_{k + 1}, I_{k + 1}, R_{k + 1}, L_{k + 1} \geq 0$. From this fact and \eqref{eq:16new} it follows the proof of the theorem.
\end{proof}
\subsection{Stability properties}
Since the discrete model \eqref{eq:15} consists of four equations, the study of its stability via the set of eigenvalues of the Jacobian is very hard, even when we use the condition $S_k + I_k + R_k + L_k = 1$ to reduce the system to a system of three equations. In the latter case the estimate of the Jacobian with the help of the Jury criterion remains complicated. As in Section 3 we shall use a generalization of Lyapunov stability theorem for getting NSFD schemes, which preserve Property $(P_3)$  of the continuous model without complex computations. For this reason consider the function
\begin{equation*}
V(I, S, L, R) = \big(I+ S - \dfrac{f}{f + g}\big)^2 + \big(L + R - \dfrac{g}{f + g}\big)^2, \qquad (I, S, L, R) \in D_4.
\end{equation*}
We shall show that the function $V$ satisfies all the conditions of  \cite[Theorem 3.2]{Iggidr} on $D_4$. Clearly, $V$ is continuous on $D_4$, where $D_4$ is defined by \eqref{eq:5}. Moreover $V(E_i^*) = 0  \; (i = 1, 2)$ and
\begin{equation*}
V(I, S, L, R) \geq 0, \quad \forall (I, S, L, R) \in D_4.
\end{equation*}
On the other hand we have
\begin{equation*}
\begin{split}
 \Delta V(I_k, S_k, L_k, R_k) &= V(I_{k + 1}, S_{k + 1}, L_{k + 1}, R_{k + 1}) - V(I_{k}, S_{k }, L_{k }, R_{k}),\\
&= \Big(I_{k + 1} + S_{k + 1} - I_k - S_k\Big)\Big(I_{k + 1} + S_{k + 1} + I_k + S_k - 2\dfrac{f}{f + g}\Big)\\
 &+ \Big(L_{k + 1} + R_{k + 1} - L_k - R_k\Big)\Big(L_{k + 1} + R_{k + 1} + L_k + R_k - 2\dfrac{g}{f + g}\Big).
\end{split}
\end{equation*}
Since Property $(P_2)$ of the model is preserved, obviously $\Delta V(I_k, S_k, L_k, R_k) \leq 0$ for any $(I_k, S_k, L_k, R_k) \in D_4$. In this case we have 
\begin{equation*}
G^* = \Big\{\big(I, S, L, R\big) \in D_4: I + S = \dfrac{f}{f + g}; \, \, L + R = \dfrac{g}{f + g}\Big\}.
\end{equation*}
In view of the fact that Property $(P_1)$ of the model is preserved all the solutions of \eqref{eq:15} are bounded.

 Thus, the conditions (1), (2) and (4) of \cite[Theorem 3.2]{Iggidr} are satisfied. It remains only to show the $G^*$- local stability of the equilibrium points  $E_i^* \, (i = 1, 2)$.\par

Notice that, for any initial conditions $\big(I_0, S_0, L_0, R_0\big)$ belonging to  $G^*$ the solution $\big(I_k, S_k, L_k, R_k\big)$ also belongs $G^*$, that is,  $I_k + S_k = f/(f + g)$ and $L_k + R_k = g/(f + g)$. Using this relation we reduce \eqref{eq:15} to two equations depending on  $I_k, L_k$
\begin{equation}\label{eq:18}
\begin{split}
&I_{k + 1} = \dfrac{\Big(1 + \varphi \beta_I \dfrac{f}{f + g} - \varphi e_I - \varphi g\Big)I_k + \varphi f L_k}{1 + \varphi \beta_I I_k},\\
 &L_{k + 1} = \dfrac{\Big(\varphi g + \varphi \beta_L \dfrac{g}{f + g}\Big)I_k + \Big(1 - \varphi f - \varphi e_L\Big)L_k}{1 + \varphi \beta_L I_k}.
\end{split}
\end{equation}
Put
\begin{equation}\label{eq:p2}
e_1^* = \big(I_1^*, L_1^* \big) = \big(0, 0\big), \qquad e_2^* = \big(I_2^*, L_2^*).
\end{equation}
It is easy to see that the equilibrium points $E_i^* \, (i = 1, 2)$ are $G^*$-locally asymptotically stable if and only if  $e_i^* (i = 1, 2)$ are locally asymptotically  stable equilibrium points of the system \eqref{eq:18}. The Jacobian of \eqref{eq:18} at $E^* = (I^*,  L^*)$ are defined by
\begin{equation}\label{eq:19}
J(E^*) = 
\begin{pmatrix}
\dfrac{1 - \varphi e_I - \varphi g + \varphi \beta_I \dfrac{f}{f + g} - \varphi^2 f \beta_I L^*}{(1 + \varphi \beta_I I^*)^2}& \quad \dfrac{\varphi f}{1 + \varphi \beta_I I^*}\\
\\
\dfrac{\varphi g + \varphi \beta_L \dfrac{g}{f + g} - \varphi \beta_L (1 - \varphi f - \varphi e_L)L^*}{(1 + \varphi \beta_L I^*)^2}& \quad \dfrac{1 - \varphi f - \varphi e_L}{1 + \varphi \beta_L I^*}
\end{pmatrix}.
\end{equation}
For convenience in the future use we restate the Jury criterion for system of two difference equations.
\begin{lemma}\label{lemma1}
The equilibrium point $E^*$ is locally asymptotically stable if and only if the Jacobian matrix $J(E^*)$ satisfies \cite[p.64, Theorem 2.10]{Linda}

\begin{enumerate}
\item $\det(J(E^*)) < 1$
\item $1 - trace(J(E^*)) + \det(J(E^*)) > 0$
\item $ 1 + trace(J(E^*)) + \det(J(E^*)) > 0$
\end{enumerate}

\end{lemma}
\textbf{(i). The first equilibrium point }
\begin{lemma}\label{lemma2}
If $c > 0$ then 
\begin{equation}\label{eq:20}
\gamma := f + e_L + e_I + g - \beta_I \dfrac{f}{f + g} > 0,
\end{equation}
where $c$ is defined by \eqref{eq:7a}.
\end{lemma}
\begin{proof}
Suppose the opposite, i.e.,  $\gamma \leq 0$.  Then
\begin{equation*}
e_I - \beta_I \dfrac{f}{f + g} \leq -f - e_L - g < -g.
\end{equation*}
From \eqref{eq:7a} we obtain
\begin{equation*}
c < (f + e_L)(-g) + g(e_L - g \beta_L \dfrac{f}{f + g}) = -fg - g \beta_L \dfrac{f}{f + g} < 0.
\end{equation*}
This is contrary to the assumption $c>0$. Thus, the lemma is proved.
\end{proof}
At the equilibrium point $e_1^*$, from \eqref{eq:19} we have
\begin{equation}\label{eq:p1}
J(e_1^*) = 
\begin{pmatrix}
1 - \varphi e_I - \varphi g + \varphi \beta_I \dfrac{f}{f + g}& \varphi f\\
\\
\varphi g + \varphi \beta_L \dfrac{g}{f + g}& 1 - \varphi f - \varphi e_L 
\end{pmatrix}.
\end{equation}

\begin{theorem}\label{theorem6}
In the case $c > 0$, if $\varphi(h)$ is a function satisfying
\begin{equation}\label{eq:21}
\varphi(h) < \min \Big\{\dfrac{\gamma}{c}, \quad \dfrac{2}{\gamma}\Big\}, \qquad \forall h > 0,
\end{equation}
then $e_1^*$ is a locally asymptotically stable equilibrium point of \eqref{eq:18}, i.e.,   $G^*$- locally asymptotically stable.
\end{theorem}

\begin{proof}
From \eqref{eq:p1} we obtain
\begin{equation*}
\det(J(e_1^*)) = 1 - \varphi \gamma + \varphi^2 c, \qquad trace(J(e_1^*)) = 2 - \varphi \gamma.
\end{equation*}
Therefore,  $\det(J(e_1^*)) < 1$ if and only if $\varphi < \gamma/c$. On the other hand 
\begin{equation*}
1 + trace(J(e_1^*)) + \det(J(e_1^*)) = 4 - 2\varphi \gamma + \varphi^2 c.
\end{equation*}
Hence, if $\varphi < 2/\gamma$ then $1 + trace(J(e_1^*)) + \det(J(e_1^*)) > 0$. Finally, we have
\begin{equation*}
1 + \det(J(e_1^*)) - trace(J(e_1^*)) = \varphi^2 c > 0.
\end{equation*}
Thus, all three conditions of Lemma
 \ref{lemma1} are satisfied. It implies that required to prove.
\end{proof}

\begin{flushleft}
\textbf{(ii). The second equilibrium point}\\
\end{flushleft}
\begin{lemma}\label{lemma3}
Consider the polynomial of degree $n$ in the variable $\varphi$ with the coefficients $a_i \in \mathbb{R}$:
\begin{equation*}
P_n(\varphi) = \sum_{i = 0}^na_i\varphi^i, \qquad \varphi > 0.
\end{equation*}
Then, if $a_0 > 0$ then there exists a number $\varphi_0 > 0$ such that
\begin{equation*}
P_n(\varphi) > 0, \qquad \forall \varphi \in (0, \varphi_0). 
\end{equation*}
\end{lemma}
\begin{proof}
The lemma is   straightforward  deduced  from the definition of the limit.
\end{proof}

\begin{lemma}\label{lamma4}
If $c < 0$ then 
\begin{equation}\label{eq:22}
I^* > \dfrac{-b + \beta_I(f + e_L)}{2a} = \dfrac{\beta_I\beta_L\dfrac{f}{f + g} - \beta_L(e_I + g)}{2\beta_I\beta_L}.
\end{equation}
\end{lemma}
\begin{proof}

From \eqref{eq:7a} we obtain
\begin{equation*}
b^2 - 4ac = -2\beta_I \beta_L c + \beta_I^2(f + e_L)^2 + (\beta_L e_I - \beta_I \beta_L \dfrac{f}{f + g})^2 + \beta_L^2(2ge_I + g^2) + 2 \beta_I \beta_L fg. 
\end{equation*}
So, if $c < 0$ we have $b^2 - 4ac > \beta_I^2(f + e_L)^2$. Having in mind this from \eqref{eq:7} we obtain \eqref{eq:22}. The lemma is proved.
\end{proof}

The Jacobian at the second equilibrium defined from \eqref{eq:19} has the form
\begin{equation}\label{eq:23a}
\det(J(e_2^*)) = \dfrac{\alpha_4 \varphi^4 + \alpha_3 \varphi^3 + \alpha_2 \varphi^2 + \alpha_1 \varphi + 1}{\big(1 + \varphi \beta_I I^*\big)^2\big(1 + \varphi \beta_L I^*\big)^2},
\end{equation}
where
\begin{equation}\label{eq:23}
\begin{split}
&\alpha_1 = \beta_L I^* - f - e_L - e_I - g + \beta_I \dfrac{f}{f + g},\\
&\alpha_2 = -\Big[(f + e_L)\beta_L I^* + (e_I + g - \beta_I \dfrac{f}{f + g})(\beta_L I^* - f - e_L)\\
&+ f\beta_IL^* + (g + \beta_L\dfrac{g}{f + g} - \beta_L L^*)f\Big],\\
&\alpha_3 = (e_I + g - \beta_I \dfrac{f}{f + g})(f + e_L)\beta_LI^* - f\beta_IL^*(\beta_L I^* - f - e_L) \\
&- (g + \beta_L\dfrac{g}{f + g} - \beta_L L^*)f\beta_II^* - (f + e_L)f\beta_LL^*, \qquad \alpha_4 = 0.
\end{split} 
\end{equation}
Its trace has the form
\begin{equation}\label{eq:24a}
trace(J(e_2^*)) = \dfrac{\gamma_4\varphi^4 + \gamma_3\varphi^3 +  \gamma_2\varphi^2 +  \gamma_1\varphi + 2}{\big(1 + \varphi \beta_I I^*\big)^2\big(1 + \varphi \beta_L I^*\big)^2},
\end{equation}
where
\begin{equation}\label{eq:24}
\begin{split}
&\gamma_1 = 3\beta_LI^* - e_I - g + \beta_I\dfrac{f}{f + g} -f - e_L + 2\beta_II^*,\\
&\gamma_2 = \beta_L^2{I^*}^2 - 2(e_I + g - \beta_I\dfrac{f}{f + g})\beta_LI^* - f\beta_IL^* + \beta_I^2{I^*}^2 + 2(\beta_L I^* - f - e_L)\beta_II^* - (f + e_L)\beta_LI^*,\\
&\gamma_3 = -(e_I + g -\beta_I\dfrac{f}{f + g})\beta_L^2{I^*}^2 - 2f\beta_I\beta_LI^*L^* +(\beta_LI^* - f - e_L)\beta_I^2{I^*}^2 - 2(f + e_L)\beta_I\beta_L{I^*}^2,\\
&\gamma_4 = -f\beta_I\beta_L^2{I^*}^2L^* - (f + e_L)\beta_L\beta_I^2{I^*}^3.
\end{split} 
\end{equation}

\begin{theorem}\label{theorem7}
If  $c < 0$ there exists a number $\varphi_0 = \varphi_0(e_I, e_L, \beta_I, \beta_L, f, g) > 0$ such that for any function $\varphi(h)$ satisfying
\begin{equation}\label{eq:25}
\varphi(h) < \varphi_0, \qquad \forall h > 0,
\end{equation}
the equilibrium point $e_2^*$ is $G^*$- asymptotically stable.
\end{theorem}

\begin{proof}
We shall use Lemma \ref{lemma1} to prove the theorem. For this reason we show that all the conditions  of the lemma are satisfied.\par
First from \eqref{eq:23a} we see that $\det(J(e_2^*)) < 1$ is equivalent to
\begin{equation}\label{eq:26}
\begin{split}
\lambda_1(\varphi) &:= \varphi^3[(\beta_I\beta_L)^2{I^*}^4 - \alpha_4] + \varphi^2[2(\beta_I + \beta_L)\beta_I\beta_L{I^*}^3 - \alpha_3]\\
 &+ \varphi[(\beta_I + \beta_L)^2{I^*}^2 + 2\beta_I\beta_L{I^*}^2 - \alpha_2] + [2(\beta_I + \beta_L)I^* - \alpha_1] > 0.
\end{split}
\end{equation}
Now we show that
\begin{equation}\label{eq:27}
A_1 := 2(\beta_I + \beta_L)I^* - \alpha_1 = \beta_LI^* + 2\beta_II^* + f + e_L + e_I + g - \beta_I\dfrac{f}{f + g} > 0.
\end{equation}
Using the estimate \eqref{eq:22} we have
\begin{equation*}
\begin{split}
2\beta_II^* > \beta_I\dfrac{f}{f + g} - (e_I + g),
\end{split}
\end{equation*}
therefore from \eqref{eq:27} we obtain $A_1 > \beta_LI^* + f + e_L > 0$. Hence, according to Lemma \ref{lemma3} there exists a number $\varphi_1 > 0$ such that for any $\varphi \in (0, \varphi_1)$ there holds $\lambda_1(\varphi) > 0$, i.e.,  $\det(J(E_2^*)) < 1$ for any $\varphi \in (0, \varphi_1)$.\\
On the other hand, combining \eqref{eq:23a} and \eqref{eq:24a} we obtain $1 - trace(J(e_2^*)) + \det(J(e_2^*)) > 0$ if and only if

\begin{equation}\label{eq:28}
\begin{split}
\lambda_2(\varphi) &:= [\beta_I^2\beta_L^2{I^*}^4 - \gamma_4 + \alpha_4]\varphi^2 + [2(\beta_I + \beta_L)\beta_I\beta_L{I^*}^3 - \gamma_3 + \alpha_3]\varphi \\
 &+ [(\beta_I + \beta_L)^2{I^*}^2 + 2\beta_I\beta_L{I^*}^2 - \gamma_2 + \alpha_2] > 0.
\end{split}
\end{equation}
We shall show that
\begin{equation*}
A_2 := (\beta_I + \beta_L)^2{I^*}^2 + 2\beta_I\beta_L{I^*}^2 - \gamma_2 + \alpha_2 > 0.
\end{equation*}
Making elementary calculations in combination with \eqref{eq:7} and \eqref{eq:7a} we obtain
\begin{equation*}
A_2 = 3a{I^*}^2 + 2bI^* + c = I^*(2aI^* + b) + (a{I^*}^2 + bI^* + c) = I^*(2aI^* + b) = I^*\sqrt{b^2 - 4ac} > 0.
\end{equation*}

Hence, according to Lemma \ref{lemma3} there exists a number $\varphi_2 > 0$ such that for any $\varphi \in (0, \varphi_2)$ there holds $\lambda_2(\varphi) > 0$, i.e., $1 - trace(J(e_2^*)) + \det(J(e_2^*)) > 0$ .\\
Finally, combining \eqref{eq:23a} and \eqref{eq:24a} we obtain $1 + trace(J(e_2^*)) + \det(J(e_2^*)) > 0$ if and only if
\begin{equation}\label{eq:28}
\begin{split}
\lambda_3(\varphi) &:= [\beta_I^2\beta_L^2{I^*}^4 + \gamma_4 + \alpha_4]\varphi^4 + [2(\beta_I + \beta_L)\beta_I\beta_L{I^*}^3 + \gamma_3 + \alpha_3]\varphi^3 \\
 &+ [(\beta_I + \beta_L)^2{I^*}^2 + 2\beta_I\beta_L{I^*}^2 + \gamma_2 + \alpha_2]\varphi^2 + [2(\beta_I + \beta_L)I^* + \gamma_1 + \alpha_1]\varphi + 4 > 0.
\end{split}
\end{equation}
Hence, according to Lemma \ref{lemma3} there exists a number  $\varphi_3 > 0$ such that for any $\varphi \in (0, \varphi_3)$ there holds  $\lambda_3(\varphi) > 0$, i.e., $1 + trace(J(e_2^*)) + \det(J(e_2^*)) > 0$ .\\
Set
\begin{equation*}
\varphi_i^* := \sup\Big\{\varphi > 0: \lambda_i(\varphi) > 0\Big\}, \quad i = 1, 2, 3; \qquad
\varphi_0 = \min_{i = 1, 2, 3}\{\varphi_i^*\}.
\end{equation*}
Then for any $\varphi < \varphi_0$ all the three conditions of Lemma \ref{lemma1} are satisfied. Thus, the theorem is proved.
\end{proof}

Summarizing the above results we obtain the following theorem of the nonstandard difference schemes preserving Properties $(P_1) - (P_3)$ of the model \eqref{eq:4}.

\begin{theorem}\label{theorem8}
Consider the model \eqref{eq:4}. Then
\begin{enumerate}
\item  In the case $c > 0$, the scheme \eqref{eq:15} preserves Properties $(P_1) - (P_3)$ of the model \eqref{eq:4} if
\begin{equation}\label{eq:29}
\varphi(h) < \varphi^* := \min\Big\{\dfrac{1}{e_I + g}, \quad \dfrac{1}{f + e_L}, \quad \dfrac{1}{f + g}, \quad \dfrac{\gamma}{c}, \quad, \dfrac{2}{\gamma}\Big\}, \qquad \forall h > 0,
\end{equation}
where $c$ and $\gamma$ are defined by \eqref{eq:7a} and \eqref{eq:20}, respectively.
\item In the case $c < 0$,  the scheme \eqref{eq:15} preserves Properties $(P_1) - (P_3)$ of the model \eqref{eq:4} if
\begin{equation}\label{eq:30}
\varphi(h) < \varphi^* : = \min\Big\{\dfrac{1}{e_I + g}, \quad \dfrac{1}{f + e_L}, \quad \dfrac{1}{f + g}, \quad\varphi_0 \Big\}, \qquad \forall h > 0,
\end{equation}
where $\varphi_0$ is defined as in Theorem \ref{theorem7}.
\end{enumerate}
\end{theorem}
\begin{remark}
There are many ways for selecting the function $\varphi$ satisfying the conditions of Theorem \ref{theorem8}, for example,  
\begin{equation*}
\varphi(h) = \dfrac{1 - e^{-\tau h}}{\tau}, \qquad \tau > \dfrac{1}{\varphi^*}.
\end{equation*}
\end{remark}

\section{Numerical simulations}
Numerical examples presented in this section show that the obtained theoretical results of the NSFD preserving the properties of the metapopulation model are valid.
\subsection{Numerical simulations for the model \eqref{eq:2}}
Very recently, in \cite{DQA} we performed numerical simulations for the model \eqref{eq:2} in particular cases. The result is that the  standard difference schemes, for example, the typical four stages Runge-Kutta do not preserve the positivity of the solution. The two stage Runge-Kutta method and the explicit Euler method have solutions oscillating near equilibrium points with increasing amplitude. Therefore, if considering the models on large time intervals, the amplitude of oscillation will be large. The numerical solutions in this case cannot preserve properties of the model \eqref{eq:2}. In general, all standard difference schemes such as  Runge-Kutta and Taylor methods preserve the properties of the continuous model for sufficiently small grid sizes.

\begin{example}{\textbf{Case $\mathcal{R}_0 < 1$}}.\\
Consider the model \eqref{eq:2} with the parameters 
\begin{equation*}
\beta = 0.8, \quad \lambda = 0.1, \quad \delta = 0.2, \quad e = 0.3.
\end{equation*}
 In this case $\mathcal{R}_0 = 0.4 < 1$. The model has two equilibrium points $P_1^* = (0.25, \, 0)$ and $P_2^*$ = (0.625, \, -0.375), among them $P_1^*$  is a globally asymptotically stable equilibrium point
 while $P_2^*$ is unstable equilibrium point.
\end{example}
For the scheme \eqref{eq:8} - \eqref{eq:9} we choose parameters $c_i$ $(i = 1, 2, 5, 6)$ satisfying Theorem \ref{Theorem5}, where
\begin{equation*}
c_1 = -17, \qquad c_2 = 18, \qquad c_5 = -1, \qquad c_6 = 2, \qquad \varphi(h) = h.
\end{equation*}
The numerical solutions obtained by the schemes in this case are depicted in Figure
 1, where each pair of blue and red curves corresponds to a solution
   $(x_k, y_k)$. In this case the properties of the model  \eqref{eq:2} are preserved.


%

\begin{example}{\textbf{Case $\mathcal{R}_0 > 1$}}.\\
Consider the model \eqref{eq:2} with the parameters 
\begin{equation*}
\beta = 2, \quad \lambda = 0.3, \quad \delta = 0.3 \quad e = 1.
\end{equation*}
 In this case $\mathcal{R}_0 = 3.75 > 1$. The model has two equilibrium points $P_1^* = (0.75, \, 0)$ and $P_2^*$ = (0.2, \, 0.55), among them $P_2^*$ is a globally asymptotically stable equilibrium point and $P_1^*$ is unstable equilibrium point.
\end{example}
For the scheme \eqref{eq:8}-\eqref{eq:9} we choose the parameters $c_i$ $(i = 1, 2, 5, 6)$ satisfying Theorem \ref{Theorem5}, where
\begin{equation*}
c_1 = -17, \qquad c_2 = 18, \qquad c_5 = -1, \qquad c_6 = 2 , \quad \varphi(h) = h.
\end{equation*}
The numerical solutions obtained by the schemes in this case are depicted in Figure 2. Clearly, the properties of the model \eqref{eq:2} are preserved.



\subsection{Numerical simulation for the model \eqref{eq:4}}
\begin{example}{\textbf{Case $c > 0$}}.\\
Consider the model \eqref{eq:4} with the parameters
\begin{equation*}
\beta_I = 0.4, \quad e_I = 0.1, \quad f = 0.25, \quad g = 0.75, \quad e_L = 1, \quad \beta_L = 0.25.
\end{equation*}
\end{example}
In this case $c =  0.7031 > 0$. The model has two equilibrium points
\begin{equation*}
P_1^* = \big(0, \, 0.25, \, 0, \, 0.75\big) \in D_4, \qquad P_2^* = \big(-1.25, \, 1.5, \, -1.25, \, 2 \big) \notin D_4,
\end{equation*}
where $P_1^*$ is a locally asymptotically stable equilibrium point, and $P_2^*$ is unstable equilibrium point.\par
The numerical solutions obtained by the four stages Runge-Kutta method and the explicit Euler method are depicted in Figures 3 and 4, respectively. Clearly, the properties of the model are not preserved.\\

The numerical solutions obtained by these differences schemes are not positive, their boundedness is destroyed. The four stages Runge-Kutta generate spurious fixed points depending on the grid sizes. The explicit Euler method gives the numerical solutions oscillating near the equilibrium points. In general, the standard difference schemes preserve the properties of the continuous model only when the grid size is sufficiently small.





We shall use NSFD schemes \eqref{eq:15} with the denominator function $\varphi(h)$ defined 
 \eqref{eq:29} in Theorem \ref{theorem8}. In this case, since $c =  0.7031$ and $\gamma = 2$ we have $\varphi^* = 0.8$. Therefore, we choose
\begin{equation*}
\varphi(h) = \dfrac{1 - e^{-2h}}{2}.
\end{equation*}
The numerical solutions of the scheme \eqref{eq:15} in this case are depicted in Figures 
 5 and 6. All the properties of the model are preserved for any grid size $h$.





\begin{example}{\textbf{Case $c < 0$}}. \\
Consider the model \eqref{eq:4} with the parameters
\begin{equation*}
\beta_I = 0.8, \quad e_I = 0.25, \quad f = 0.2, \quad g = 0.75, \quad e_L = 0.1, \quad \beta_L = 2.
\end{equation*}
\end{example}
In this case we have $c = -0.2163 < 0$. The model has two equilibrium points

\begin{equation*}
P_1^* = \big(0, \, 0.2105, \, 0, \, 0.7895\big) \in D_4, \qquad P_2^* = \big(0.1045, \, 0.1060, \,  0.4781, \, 0.3114 \big) \in D_4,
\end{equation*}
where $P_1^*$ is an unstable equilibrium point and  $P_2^*$ is a locally asymptotically stable equilibrium point.\par
As in the previous examples, the standard difference schemes such as Runge-Kutta or Taylor methods  cannot preserve the properties of the continuous model. We shall use NSFD schemes \eqref{eq:15} with the denominator function $\varphi(h)$ defined by \eqref{eq:30} in Theorem \ref{theorem8}. In this case we have
\begin{equation*}
\begin{split}
\lambda_1(\varphi) &= 3.0509 \times 10^{-4} \varphi^3 +  0.0314 \varphi ^2 + 0.4586 \varphi + 1.5077,\\
\lambda_2(\varphi) &=  0.0041  \varphi^2 +  0.0684  \varphi + 0.2338,\\
\lambda_3(\varphi) &=  -0.0035 \varphi^4  - 0.0904 \varphi^3 - 0.6688 \varphi^2 - 0.6750 \varphi + 4.
\end{split}
\end{equation*}
The graphs of the functions $\lambda_i \; (i = 1, 2, 3)$ are given in Figure 7. We can choose the number $\varphi_0$ in Theorem \ref{theorem8} equal to $1. 5$. Using Theorem \ref{theorem8} we choose $$\varphi(h) = \dfrac{1 - e^{-1.1 h}}{1.1}.$$ The numerical solutions obtained in this case are depicted in Figures 8 and 9. From the figures we see that all the properties of the model are preserved.

%
%





\section{Conclusion}
In this paper we have constructed the NSFD schemes preserving the properties of two metapopulation models. An essential result in the selection of the parameters of the schemes is that we used a generalization of Lyapunov stability theorem to ensure the stability properties of the the schemes. This approach is much simpler and more effective than the approach of studying Jacobians of the discrete systems which was used by ourselves in a previous work  for the first metapopulation model \cite{DQA}. In the future we shall develop this approach for constructing discrete models preserving the properties of continuous models for other problems.

\section*{Acknowledgments}
This work is supported by Vietnam National Foundation for Science and Technology
Development (NAFOSTED) under the grant  number 102.01-2014.20. \\


\begin{thebibliography}{99}
\bibitem{Linda} 
L. J. S. Allen, {\em An Introduction to Mathematical Biology}, Prentice Hall, New Jersey, $2007$.
\bibitem{Amarasekare}
P. Amarasekare, H. Possingham, {\em Patch Dynamics and Metapopulation Theory: the Case of  Successional Species},  Journal of Theoretical Biology, 209 (2001) pp. 333-344.
\bibitem{AL1}
R. Anguelov, J. M. -S Lubuma, {\em Nonstandard finite difference method by nonlocal approximations}, Mathematics and Computers in Simulation, {61} $(2003)$,  pp. 465-475.
\bibitem{AL2} 
R. Anguelov, Y. Dumont, J.M.-S. Lubuma, M. Shillor, {\em Dynamically consistent nonstandard finite difference schemes for epidemiological models}, Journal of Computational and Applied Mathematics, {255} (2014), pp. 161-182.
\bibitem{Basak}
S. C. Basak, G. Restrepo, J. L. Villaveces, {\em Advances in Mathematical Chemistry and Applications}, Bentham Science Publishers, 2016.
\bibitem{Barrante} 
J. R. Barrante, 	{\em Applied Mathematics for Physical Chemistry}, Prentice Hall, 2003.
\bibitem{Brauer} 
F. Brauer, C. Castillo - Chavez, {\em Mathematical Models in Population Biology and Epidemiology}, Springer New York, $(2001)$.
\bibitem{DQA}
Q. A. Dang, M. T. Hoang, {\em Dynamically consistent discrete metapopulation model}, Journal of Difference Equations and Applications, 2016, DOI: 10.1080/10236198.2016.1197213.
\bibitem{DK1} 
T. D. Dimitrov, H. V. Kojouharov, {\em Stability-Preserving Finite-Difference Methods For General Multi-Dimensional Autonomous Dynamical Systems}, {International Journal Of Numerical Analysis And Modeling}, {4} (2) (2007), pp. 280-290.
\bibitem{DK2} 
T. D. Dimitrov, H. V. Kojouharov, {\em Nonstandard finite difference schemes for general two - dimensional autonomous dynamical systems}, Applied Mathematics Letters, {18} (2005), pp. 769-774.

\bibitem{DK3}
 T. D. Dimitrov, H. V. Kojouharov, {\em Positive and elementary stable nonstandard numerical methods with applications to predator - prey models}, Journal of Computational and Applied Mathematics, {189} (2006), pp.  98-108.
\bibitem{DK4}D. T. Dimitrov \& H. V. Kojouharov, {\em Dynamically consistent numerical methods for general productive–destructive systems}, Journal of Difference Equations and Applications, 17(12) (2011), pp. 1721-1736.
\bibitem{DK5}
 T. D. Dimitrov, H. V. Kojouharov, {\em Nonstandard finite-difference methods for predator-prey models with general functional response}, Mathematics and Computers in Simulation, {78} (2008), pp. 1-11.

\bibitem{E1} 
M. Ehrhardt, R.E. Mickens,  {\em A Nonstandard Finite Difference Scheme for Convection-Diffusion Equations having Constant Coefficients}, Appl. Math. Comput. {219} (2005), pp. 6591-6604.
\bibitem{Elaydi}
S. Elaydi, {\em An Introduction to Difference Equations}, Springer Science+Business Media, Inc, 2005. 

\bibitem{Garba}S.M. Garba, A.B. Gumel, J.M.-S. Lubuma, {\it Dynamically-consistent non-standard finite difference method for an epidemic model}, Mathematical and Computer Modelling {53} (2011), pp. 131-150.

\bibitem{Henner} 
V. Henner, T. Belozerova, K. Forinash, {\em Mathematical Methods in Physics: Partial Differential Equations, Fourier Series, and Special Functions}, A K Peters, Wellesley, 2013.

\bibitem{Iggidr}
A. Iggidr, M. Bensoubaya, {\em New Results on the Stability of Discrete-Time Systems and Applications to Control Problems}, Journal of Mathematical Analysis and Applications, {219} (1998),  pp. 392-414.

\bibitem{Kaczor}
W. J. Kaczor, M. T. Nowak, {\em Problems in Mathematical Analysis I: Real Numbers, Sequences and Series}, American Mathematical Society, 2000.
\bibitem{Keymer}
J.E. Keymer, P.A. Marquet, J. X. Velasco-Hernandez, . S. A. Levin, {\em Extinction thresholds and metapopulation persistence in dynamic landscapes}, The American Naturalist, Vol. 156 (2000), No. 5, pp. 478-494.

\bibitem{Keshet} 
L. Edelstein-Keshet, {\em Mathematical Models in Biology}, Society for Industrial and Applied Mathematics, Philadelphia, 1998.
\bibitem{K1}  
H. V. Kojouharov, D. T. Wood, {\em A class of nonstandard numerical methods for autonomous dynamical systems}, Applied Mathematics Letters, {50} (2015), pp. 78-82.
\bibitem{K2} H. Kojouharov , B. Welfert, {\em A nonstandard Euler schemes for $y'' + g(y)y' + f(y)y = 0$}, Journal Computational and Applied Mathematics, {151} (2003), pp. 335-353.
\bibitem{Leveque}
R. Leveque, 	{\em Finite Difference Methods for Ordinary and Partial Differential Equations}, Classics in Applied Mathematics, Society for Industrial and Applied Mathematics, 2007.
\bibitem{Marchuk1} 
G. I. Marchuk, {\em Mathematical Modelling of Immune Response in Infectious Diseases}, Springer, 1997.
\bibitem{Marchuk2}
G.I. Marchuk, {\em Mathematical Modelling in Environmental Problems, in: Studies in Mathematics and its Applications}, vol. 16, North-Holland, 1986.


\bibitem{Mickens0}
R. E. Mickens, {\em “Exact solutions to a finite-difference model of a nonlinear reaction-advection equation: Implications for numerical analysis”}, Numerical Methods for Partial Differential Equations, 5 (1989), pp. 313–325.

\bibitem{Mickens1} 
R. E. Mickens, {\em Nonstandard Finite Difference Models of Differential Equations}, World Scientific, Singapore, (1994).

\bibitem{Mickens2} 
R. E. Mickens, {\em Applications  of Nonstandard Finite Difference Schemes}, World Scientific, Singapore, (2000).
\bibitem{Mickens3} 
R. E. Mickens, {\em Advances in the Applications of Nonstandard Finite Difference Schemes}, World Scientific, Singapore, New Jersey, 2005.
\bibitem{Mickens4} 
R. E. Mickens,  {\em Nonstandard Finite Difference Schemes for Differential Equations}, Journal of Difference Equations and Applications, {8}(9) (2005), pp. 823-847.
\bibitem{Mickens5}  R. E. Mickens, { \em Dynamic consistency: a fundamental principle for constructing nonstandard finite difference schemes for differential equations}, Journal of Difference Equations and Applications, 11(7) (2005), pp. 645-653.
\bibitem{Partidar}
K. C. partiadar, {\em Nonstandard finite difference methods: recent trends and further developments}, Journal of Difference Equations and Applications, 2016, DOI:10.1080/10236198.2016.1144748.
\bibitem{Petzold}
L. R. Petzold, {\em Computer methods for Ordinary Differential Equations and Differential-Algebraic Equations}, Society for Industrial and Applied Mathematics, 1998.
\bibitem{Robinson}
J. C. Robinson, J. L. Rodrigo, 	{\em Partial Differential Equations and Fluid Mechanics, London Mathematical Society Lecture Note Series}, Cambridge University Press, 2009.



\bibitem{Roeger1} 
L.-I. W. Roeger,  {\em General nonstandard finite-difference schemes for differential equations with three fixed-points}, Computers and Mathematics with Applications, {57} (2009), pp. 379-383.
\bibitem{Roeger2}
  L.-I. W. Roeger, {\em Nonstandard finite difference schemes for differential equations with $n + 1$ distinct fixed-points}, Journal of Difference Equations and Applications, {15} (2009), pp. 133-151.
\bibitem{Roeger3} 
L.-I. W. Roeger, {\em Dynamically Consistent Discrete - Time Lotka-Volterra Competition Models}, Discrete And Continuous Dynamical Systems,  Supplement 2009, pp. 650-658.
\bibitem{Roeger4} 
L.-I. W. Roeger, {\em Periodic solutions preserved by nonstandard finite-difference schemes for Lotka-Volterra system: a different approach}, Journal of Difference Equations and Applications, {14} (2008), pp. 481-493.
\bibitem{Roeger5} 
L.-I. W. Roeger, {\em Nonstandard finite-difference schemes for the Lotka-Volterra systems: generalization of Mickens's method}, Journal of Difference Equations and Applications, {12} (2006), pp. 937-948.
\bibitem{Roeger6} L.-I.W. Roeger, {\it Dynamically consistent discrete-time SI and SIS epidemic models}, Discrete and Continuous Dynamical Systems, Supplement (2013), pp. 653-662.

\bibitem{Samarskii}
A. A. Samarskii, {\em The theory of difference schemes, Marcel Dekker}, New York, 2001.
\bibitem{Smith} 
H. L. Smith, P. Waltman, {\em The Theory of the Chemostat}, Cambridge Univ. Press, Cambridge, U.K. 1995.
\bibitem{Stakgold} 
I. Stakgold, D. D. Joseph, D. H. Sattinger, {\em Nonlinear Problems in the Physical Sciences and Biology}, Lecture Notes in Mathematics, Springer, 1973.
\bibitem{Strikwerda}
J. Strikwerda, {\em Finite Difference Schemes and Partial Differential Equations}, SIAM Filadelphia, 2007.

\bibitem{Wood}
D. T. Wood, D. T. Dimitrov,  Hristo V. Kojouharov, {\em A nonstandard finite difference method for n-dimensional productive-destructive systems}, Journal of Difference Equations and Applications, 2016, DOI: 10.1080/10236198.2014.997228.

\end{thebibliography}
\end{document}